\documentclass[12pt,a4paper,tbtags,reqno]{amsart}
\usepackage{amssymb}
\usepackage[space,nosort,noadjust]{cite}
\usepackage[width=160mm]{geometry}

\allowdisplaybreaks

\newtheorem{thm}{Theorem}

\newtheorem{cor}{Corollary}

\theoremstyle{definition}
\newtheorem{utv}{Proposition}
\newtheorem*{dfn}{Definition}

\theoremstyle{remark}
\newtheorem{rem}{Remark}

\numberwithin{equation}{section}

\newcommand{\vG}{\varGamma}

\newcommand{\pwp}{\wp^{\prime}}
\newcommand{\wta}{\widetilde{a}}
\newcommand{\wtc}{\widetilde{c}}
\newcommand{\wtd}{\widetilde{d}}

\newcommand{\wtI}{\widetilde{I}}

\newcommand{\wtu}{\widetilde{u}}
\newcommand{\wtE}{\widetilde{E}}
\newcommand{\wtG}{\widetilde{\vG}}
\newcommand{\whPsi}{\widehat{\Psi}}

\newcommand{\whI}{\widehat{I}}

\newcommand{\Bz}{\mathbb{Z}}

\newcommand{\Bc}{\mathbb{C}}
\newcommand{\ii}{\mathrm{i}}
\newcommand{\const}{\mathrm{const}}
\newcommand{\Ll}{\mathcal{L}}

\DeclareMathOperator*{\Res}{Res}
\def\nodot#1#2{}

\begin{document}

\title{Finite-gap solutions of the Fuchsian equations}

\author[A. O. Smirnov]{Alexander O. Smirnov}
\address{Department of Mathematics
St.-Petersburg State University of Aerospace Instrumentation
Bolshaya Morskaya str. 67, St.-Petersburg 190000, Russia}
\email{alsmir@peterlink.ru}
%

\begin{abstract}
We find a new class of the Fuchsian equations, which have an
algebraic geometric solutions with the parameter belonging to a
hyperelliptic curve. Methods of calculating the algebraic genus of
the curve, and its branching points, are suggested. Numerous
examples are given.
\end{abstract}

\maketitle

\section*{Introduction}

Integrability of the Heun equation with half-odd characteristic
exponents
\begin{equation}
\dfrac{d^2y}{dz^2}+P(z)\dfrac{dy}{dz}+Q(z)y=0, \label{eq:main}
\end{equation}
where
\begin{align}
P(z)&=\frac12\left(\frac{1-2m_1}z+\frac{1-2m_2}{z-1}
+\frac{1-2m_3}{z-a}\right), \label{eq:P:heun}\\
Q(z)&=\frac{N(N-2m_0-1)z+\lambda}{4z(z-1)(z-a)},
\label{eq:Q:heun}\\
N&=m_0+m_1+m_2+m_3, \qquad m_i \in\Bz_{\geqslant0},\qquad \lambda,
z\in\Bc, \label{eq:N:heun}
\end{align}
%
was probably discovered by Darboux more then 100 years ago
\cite{D882}. But only recently so-called finite-gap solutions
\begin{equation} \label{sol:heun}
Y_{1,2}(\mathbf{m};\lambda;z)=
\sqrt{\Psi_{g,N}(\lambda,z)}\exp\left(\pm\frac{\ii\nu(\lambda)}2
\int\frac{z^{m_1}(z-1)^{m_2}(z-a)^{m_3}\,d
z}{\Psi_{g,N}(\lambda,z) \sqrt{z(z-1)(z-a)}}\right)
\end{equation}
of this Heun equation were wrote out and analyzed \cite{Sm2002,
Tak}.
Here $\ii^2=-1$,
\begin{equation*}
\vG:\quad \nu^2=\prod_{j=1}^{2g+1}(\lambda-\lambda_j),
\end{equation*}
 $\Psi_{g,N}(\lambda,z)$ is some
polynomial of the degree $N$ in $z$ and of the degree $g$ in
$\lambda$.

In works  \cite{D882,Sm2002, Tak} the connection between Heun
equation \eqref{eq:main}--\eqref{eq:N:heun} and Treibich-Verdier
equation
\begin{equation}
\psi_{xx}-u(x)\psi=E\psi, \label{eq:shr}
\end{equation}
\begin{equation}
u(x)=m_0(m_0+1)\wp(x)+\sum_{i=1}^3m_i (m_i+1)\wp(x-\omega_i),
\label{pot:tv}
\end{equation}
was used.
%
Here $\wp(x)$ is the Weierstrass function \cite{Akh},
\begin{gather}
[\pwp(x)]^2=4\wp^3(x)-g_2\wp(x)-g_3=4\prod_{j=1}^3(\wp(x)-e_j),
\label{eq:torus}\\
\wp(\omega_i)=e_i, \qquad \wp(x-2\omega_i)\equiv \wp(x). \notag
\end{gather}

The equation \eqref{eq:shr}, \eqref{pot:tv} is generalization of
well-known Lam\'e equation
\begin{equation}
\psi_{xx}-N(N+1)\wp(x)\psi=E\psi. \label{eq:lame}
\end{equation}

At beginning of 90-th author \cite{Smr89a,Smr94} investigated
Shr\"odinger operator with finite-gap elliptic potentials and
proposed the next finite-gap elliptic generalization of
potential\eqref{pot:tv}:

%
\begin{multline}
u(x)=m_0(m_0+1)\wp(x)+\sum_{i=1}^3m_i (m_i+1)\wp(x-\omega_i)+\\
+\sum_{k=1}^M
n_k(n_k+1)\left\{\wp(x-\delta_k)+\wp(x+\delta_k)\right\}
\label{pot:gen}
\end{multline}
%
Treibich later \cite{Tr99, Тр2001} proved that for $M=1$, $n_1=1$
and for any system of numbers $m_i\in\Bz_{\geqslant 0}$ there
exists a point $\delta_1$ such that potential \eqref{pot:gen} is
finite-gap potential of Shr\"odinger operator \eqref{eq:shr}. On
the other hand the formula\footnote{prime after sign of summation
as ever indicate what term with $i^2+j^2=0$ is omitted.}
\begin{multline} \label{eq:multi}
\wp\left(x\left| \frac{2\omega}k,\frac{2\omega'}l \right)\right.=
\sum_{i=0}^{l-1}\sum_{j=0}^{k-1}
\wp\left.\left(x+\dfrac{2j\omega}{k}+ \dfrac{2i\omega'}{l}\right|
2\omega,2\omega'\right)-\\
-\sum\nolimits' \wp\left.\left(\dfrac{2j\omega}{k}+
\dfrac{2i\omega'}{l}\right| 2\omega,2\omega'\right),
\end{multline}
%
where
\begin{equation*}
\omega=\omega_1,\quad  \omega'=\omega_3,\quad
\omega_2=\omega+\omega',
\end{equation*}
and examples in works \cite{Smr89a,Smr94} show us what Treibich
cases do not exhaust all the set of even elliptic finite-gap
potentials (see for example remark~\ref{rem3}).

Since parallel using of properties of algebraic polynomials and of
elliptic meromorphic functions make possible a big progress in
study Heun and Treibich-Verdier equations (in particular simple
methods of obtaining spectral curves and monodromy matrix
\cite{Sm2002, Tak}), we decide to apply this method to analysis of
Shr\"odinger equation \eqref{eq:shr} with
potentials\eqref{pot:gen} and of appropriate Fuchsian equations.

In present work we introduce concept of `finite-gap' Fuchsian
equation, and prove necessary and sufficient conditions of
`finite-gapness' of Fuchsian equation. Also we ge\-ne\-ra\-li\-ze
the formula \eqref{sol:heun} for the case `finite-gap' Fuchsian
equation with five and more singular points, and derive the
equation of appropriate spectral curve. Algebraic genus of
spectral curve is bounded and calculated. The equation with five
singular point is examined in detail. For equation with five
singular points ($M=1$, $n_1=1$) the condition for position of
fifth point is found.

The needed facts from the theory of finite-gap elliptic potentials
for the Schr\"o\-dinger operator are collected without proof in
the first paragraph.

In the Appendix we give several simplest solutions of the
`finite-gap' Fuchsian equation and appropriate finite-gap
potentials and their spectral curves.

The author thanks A.~Treibich, V.B.~Matveev, V.B.~Kuznetsov and
V.Z.~Enolski\u\i\ for useful discussions.


\section{Schr\"odinger operator with finite-gap elliptic potential}

\begin{utv}[\cite{Nov74,ZMNP,KD,Mat76}]
Any $g$-gap potential $u(x)$ of the Schr\"odinger operator
\eqref{eq:shr} satisfies the Novikov equation:
\begin{equation}
J_g+\sum_{m=1}^g c_m J_{g-m}=d, \label{eq:novikov}
\end{equation}
%
or, which is the same, the stationary `higher' Korteweg-de Vries
(KdV) equation:
\begin{equation*}
\partial_x\left(J_g+\sum_{m=1}^g c_m J_{g-m}\right)=0.
\end{equation*}
Here $c_m$, $d$ are constants and the functions $J_m$ are the
flows of the `higher' KdV equations
%
\begin{equation*}
\partial_{t_m}u=\partial_x J_m.
\end{equation*}
The expressions for the flows $J_m$ are found from the relations

%
\begin{subequations} \label{kdv:int:J}
\begin{align}
&L\psi=\psi_{xx}-4u\psi+2u_x\int_{x}^{\infty}\psi(\tau)d\tau, \\
&(J_n)_x=L^n(u_x).
\end{align}
\end{subequations}
%
where $u(x)$ is a potential decreasing fast at $\infty$.
\end{utv}

In particular,
\begin{align*}
&J_0=u, \qquad  J_1=u_{xx}-3u^2, \qquad
J_2=u_{xxxx}-10u_{xx}u-5u^2_x+10u^3, \\
&J_3=u_{xxxxxx}-14uu_{xxxx}-28u_xu_{xxx}-21u_{xx}^2
+70u^2u_{xx}+70uu^2_x- 35u^4.
\end{align*}

\begin{rem}
In the case of decreasing fast at $\infty$ potential $u(x,t)$ the
variables
\begin{equation*}
C_j(t)=\int_{-\infty}^{\infty} J_j(x,t)\,d x
\end{equation*}
%
constitute an infinite set of integrals of motion for the KdV
equation
\begin{equation*}
\partial_t u=\partial_x J_1.
\end{equation*}
\end{rem}

\begin{utv}[\cite{DMN76,Mat76,ZMNP}]
The function
\begin{equation}
\whPsi(x,E)=E^g+\sum_{j=1}^g \gamma_j(x) E^{g-j},
\label{kdv:prod:psi}
\end{equation}
%
where
\begin{equation}
\gamma_j(x)=-\frac2{4^j}\left(J_{j-1} +\sum_{m=1}^{j-1}c_m
J_{j-m-1}-\frac{c_j}2\right), \label{kdv:gamma}
\end{equation}
obeys the equation
%
\begin{equation*}
\whPsi_{xxx}=4(u+E)\whPsi_x+2u_x\whPsi,
\end{equation*}
%
the solutions of which are the products of any two solutions of
the equation \eqref{eq:shr}. Here $c_j$, $J_j$ are the same as in
\eqref{eq:novikov} and $u(x)$ is a $g$-gap potential of the
Schr\"odinger operator \eqref{eq:shr}.
\end{utv}

\begin{dfn}[\cite{GW94b,GW94c}]
Let $u(x)$ be an elliptic function. If the general solution of
equation  \eqref{eq:shr} is meromorphic for each complex number
$E$ then $u(x)$ is called a Picard potential.
\end{dfn}

\begin{utv}[\cite{GW94b,GW94c}] \label{utv:picard}
Function $u(x)$ is an elliptic finite-gap potential if and only if
$u(x)$ is a Picard potential.
\end{utv}

\section{
Even elliptic potentials of Shr\"odinger operator and Fuchsian
equation}

Let us consider the Shr\"odinger equation \eqref{eq:shr} with even
elliptic potential \eqref{pot:gen}. It is not difficult to check
that with the help of appropriate change of variable this equation
can be transformed into partial case of Fuchsian equation with
$M+4$ singular points.

For $M=0$ equation of change variable \cite{Sm2002}\footnote{In
\cite{D882}, \cite{Tak} another changes are used}
\begin{equation} \label{eq:x2z:heun}
\begin{gathered}
\psi(x)=y(z)z^{-m_1/2}(z-1)^{-m_2/2}(z-a)^{-m_3/2},\quad
\wp(x)=e_1+(e_2-e_1)z,\\
e_2=\frac{a-2}{a+1}e_1,\qquad e_3=\frac{1-2a}{a+1}e_1,\quad
E=(e_1-e_2)\lambda+\const.
\end{gathered}
\end{equation}
%
transform equation \eqref{eq:shr}, \eqref{pot:tv} into
Heun equation \eqref{eq:main}-\eqref{eq:N:heun}.

For $M=1$ following equation of change variable
\begin{equation} \label{eq:x2z:five}
\begin{gathered}
\psi(x)=y(z)z^{-m_1/2}(z-1)^{-m_2/2}(z-a)^{-m_3/2}(z-b)^{-n_1},\\
\wp(x)=e_1+(e_2-e_1)z,\quad
E=(e_1-e_2)\lambda+\const,\\
e_2=\frac{a-2}{a+1}e_1,\qquad e_3=\frac{1-2a}{a+1}e_1,\quad
\wp(\delta_1)=\frac{a+1-3b}{a+1}e_1,
\end{gathered}
\end{equation}
%
transform equation \eqref{eq:shr}, \eqref{pot:gen} into Fuchsian
equation \eqref{eq:main} with five singular points, where
\begin{align}
P(z)&=\frac12\left(\frac{1-2m_1}z+\frac{1-2m_2}{z-1}
+\frac{1-2m_3}{z-a}\right)-\frac{2n_1}{z-b}, \label{eq:P:five}\\
Q(z)&=\frac{N(N-2m_0-1)z+\lambda+2n_1\rho(z-b)^{-1}}{4z(z-1)(z-a)},
\label{eq:Q:five}\\
\rho&=2(b-1)(b-a)m_1+2b(b-a)m_2+2b(b-1)m_3+(3b^2-2(a+1)b+a)n_1,
\label{eq:rho:five}\\
N&=m_0+m_1+m_2+m_3+2n_1, \qquad m_i,n_1 \in\Bz_{\geqslant0},\qquad
\lambda, z\in\Bc. \label{eq:N:five}
\end{align}

It is not difficult to generalize this change of variables for
greater number of singular points
\begin{equation} \label{eq:x2z:six}
\begin{gathered}
\psi(x)=y(z)z^{-m_1/2}(z-1)^{-m_2/2}(z-a)^{-m_3/2}\prod_{k=1}^M(z-b_k)^{-n_k},\\
\wp(x)=e_1+(e_2-e_1)z,\quad E=(e_1-e_2)\lambda+\const\\
e_2=\frac{a-2}{a+1}e_1,\qquad e_3=\frac{1-2a}{a+1}e_1,\quad
\wp(\delta_k)=\frac{a+1-3b_k}{a+1}e_1.
\end{gathered}
\end{equation}
%
Coefficients of corresponding Fuchsian equation \eqref{eq:main}
with $M+4$ singular points has a next form
\begin{align}
P(z)=&\frac12\left(\frac{1-2m_1}z+\frac{1-2m_2}{z-1}
+\frac{1-2m_3}{z-a}\right)-2\sum_{k=1}^M \frac{n_k}{z-b_k}, \label{eq:P:six}\\
Q(z)=&\frac{N(N-2m_0-1)z+\lambda+R_M(z)}{4z(z-1)(z-a)},\quad
R_M(z)=2\sum_{k=1}^M \frac{n_k\rho_k}{z-b_k}, \label{eq:Q:six}\\
\rho_k=&2(b_k-1)(b_k-a)m_1+2b_k(b_k-a)m_2+2b_k(b_k-1)m_3+ \notag\\
&+(3b_k^2-2(a+1)b_k+a)n_k
+4\sum_{j\ne k}\dfrac{b_k(b_k-1)(b_k-a)}{b_k-b_j}n_j, \label{eq:rho:six}\\
N=&\sum_{j=0}^3 m_j+2\sum_{k=1}^M n_k, \qquad m_j,n_k
\in\Bz_{\geqslant0},\qquad \lambda, z\in\Bc.  \label{eq:N:six}
\end{align}

\section{Finite-gap Fuchsian equations}

Let us introduce the concept of `finite-gap' solutions of the
Fuchsian equation.

\begin{dfn}
The Fuchsian equation \eqref{eq:main},
\eqref{eq:P:six}--\eqref{eq:N:six} will be called
\emph{`finite-gap'} if it can be obtained from Shr\"odinger
equation with finite-gap potential \eqref{pot:gen} via the change
of variables \eqref{eq:x2z:six}. Solutions of `finite-gap'
Fuchsian equations will be called \emph{`finite-gap' solutions}.
\end{dfn}

\begin{thm} \label{thm1}
The Fuchsian equation \eqref{eq:main},
\eqref{eq:P:six}--\eqref{eq:N:six} is `finite-gap' if and only if
points $z=b_k$ are false singular points.
\end{thm}

\begin{proof}
Let us calculate characteristic exponents at singular points of
equation
\begin{align*}
&\rho_1(0)=\dfrac{1}{2}+m_1,\qquad \rho_2(0)=0,\\
&\rho_1(1)=\dfrac{1}{2}+m_2,\qquad \rho_2(1)=0,\\
&\rho_1(a)=\dfrac{1}{2}+m_3,\qquad \rho_2(a)=0,\\
&\rho_1(\infty)=-\dfrac{N}{2},\qquad \rho_2(\infty)=-\dfrac{N-2m_0-1}{2},\\
&\rho_1(b_k)=2n_k+1,\qquad \rho_2(b_k)=0.
\end{align*}

Assume that Fuchsian equation is `finite-gap' but points $z=b_k$
are not false singular points. The a general solution will have at
those points logarithmic singularities (see for example
\cite{Ince, ODE}). Hence a general solution of Shr\"odinger
equation \eqref{eq:shr} with finite-gap potential \eqref{pot:gen}
will have logarithmic singularities at points $ x=\pm\delta_k$.
That is in contradiction with proposition~\ref{utv:picard}.

If we now assume that points $ z=b_k $ are false singular points
of equation \eqref{eq:main}, \eqref{eq:P:six}--\eqref{eq:N:six}
then from properties of general solution of this equation and from
the equation of the change of variable \eqref{eq:x2z:five} it
follows that singularities of general solution  $\psi(x)$ of
equation \eqref{eq:shr}, \eqref{pot:gen} are only poles:
\begin{align*}
&\psi(x)=c_1(x-\omega_1)^{-m_1}+o\left((x-\omega_1)^{-m_1}\right),\qquad
x\to\omega_1;\\
&\psi(x)=c_2(x-\omega_2)^{-m_2}+o\left((x-\omega_2)^{-m_2}\right),\qquad
x\to\omega_2;\\
&\psi(x)=c_3(x-\omega_3)^{-m_3}+o\left((x-\omega_3)^{-m_3}\right),\qquad
x\to\omega_3;\\
&\psi(x)=d_k(x\mp\delta_k)^{-n_k}+o\left((x\mp\delta_k)^{-n_k}\right),\qquad
x\to\pm\delta_k.
\end{align*}
%
I.e. a general solution $\psi(x)$ is a meromorphic function and
potential \eqref{pot:gen} is Picard potential. Therefore potential
\eqref{pot:gen} is finite-gap and Fuchsian equation is
`finite-gap' equation.
\end{proof}

Now let us follow \cite{Sm2002} and give the definition of
Novikov's equation for the Fuchsian equation.

\begin{dfn} By {\it Novikov's equation
of order $g$ for the Fuchsian equation} we will call the following
equality:
%
\begin{equation}
I_g+\sum_{j=1}^g \wtc_j I_{g-j}=\wtd, \label{eq:nov1}
\end{equation}
%
where $\wtc_j$, $\wtd$ are some constants,
\begin{subequations}\label{six:I}
\begin{equation}
I_{j+1}=\Ll(I_j),
\end{equation}
\begin{equation}
\Ll(f)=z(z-1)(z-a)\frac{d^2f}{dz^2}
+\frac{3z^2-2(a+1)z+a}2\cdot\frac{d f}{d z}- \int\left(4I_0\frac{d
f}{d z}+2f\frac{dI_0}{d z}\right)d z,
\end{equation}
\end{subequations}
\begin{multline} \label{six:I0}
\qquad I_0(z)=\frac{m_0(m_0+1)}4z +\frac{m_1(m_1+1)}4\cdot\frac a z+{}\\
{}+\frac{m_2(m_2+1)}4\cdot\frac{z-a}{z-1}
+\frac{m_3(m_3+1)}4\cdot\frac{a(z-1)}{z-a}+\\
+\sum_{k=1}^M n_k(n_k+1)\left(\dfrac{b_k(b_k-1)(b_k-a)}{(z-b_k)^2}
+\dfrac{3b_k^2-2(a+1)b_k+a}{2(z-b_k)}\right).
\end{multline}
\end{dfn}

\begin{thm}
The function $I_0$ \eqref{six:I0} satisfies Novikov's equation of
order $g$ \eqref{eq:nov1}, \eqref{six:I} for the Fuchsian equation
if and only if the potential $u(x)$ \eqref{pot:gen} satisfies
Novikov's equation \eqref{eq:novikov}, \eqref{kdv:int:J} of the
same order.\footnote{cf. \cite{Sm2002}}
\end{thm}

\begin{rem} \label{rem1}
The constant of integration in the definition of the function
$I_n$ \eqref{six:I} is not fixed. Because of that, all the
functions $I_n$ are defined modulo a linear combination of lower
order functions $I_k$. However, it is easily seen that the
property of `finite-gapness' of the Fuchsian equation does not
depend on the concrete values of integration constants in the
definition of $I_n$. The latter affect only the values of the
constants $\wtc_m$ and $\wtd$ in the equation \eqref{eq:nov1}.

\end{rem}

\begin{cor}
The Fuchsian equation \eqref{eq:main},
\eqref{eq:P:six}--\eqref{eq:N:six} is `finite-gap' if and only if
the function $I_0$ \eqref{six:I0} satisfies Novikov's equation of
order $g$ \eqref{eq:nov1}, \eqref{six:I} for the Fuchsian
equation.
\end{cor}

\begin{thm} \label{thm3}
For any $m_i\in\Bz_{\geqslant 0}$ there exists number $b$ such
that the Fuchsian equation with five singular points
\eqref{eq:main}, \eqref{eq:P:five}--\eqref{eq:N:five} for $n_1=1$
is `finite-gap'.
\end{thm}

\begin{proof}
From the properties of the flows $J_n$ \eqref{kdv:int:J} (see,
e.g., \cite{KD}), from the properties of elliptic functions
\cite{Akh} and from the equation of the change of variable
\eqref{eq:x2z:five} it follows that all functions $I_n$
\eqref{six:I} are rational functions of the variable $z$ (i.e.
these functions do not have logarithmic singularities).

From equation \eqref{six:I0} it follows what $I_0$ has at the
point $z=b$ a pole of second order at the point $z=b$
\begin{equation} \label{asym:I0}
I_0=n_1(n_1+1)\left(\dfrac{b(b-1)(b-a)}{(z-b)^2}
+\dfrac{3b^2-2(a+1)b+a}{2(z-b)}\right)+O(1),\quad z\to b.
\end{equation}

If we now assume that the function $I_j(z)$ has at the point $z=b$
a pole of order $2\alpha\leqslant 2n_1$:
\begin{equation*}
I_j(z)=\frac A{(z-b)^{2\alpha}}+\frac
B{(z-b)^{2\alpha-1}}+O\left((z-b)^{2-2\alpha}\right), \qquad z\to
b,
\end{equation*}
%
then we obtain that the function $I_{j+1}(z)$,
\begin{multline*}
I_{j+1}=\dfrac{2(2\alpha+1)(\alpha+n_1+1)(\alpha-n_1)}{\alpha+1}\cdot
\dfrac{b(b-1)(b-a)A}{(z-b)^{2\alpha+2}}+{}\\
{}+\dfrac{(4\alpha+1)(2\alpha^2+\alpha-n_1-n_1^2)}{2\alpha+1}\cdot
\frac{(3b^2-2(a+1)b+a)A}{(z-b)^{2\alpha+1}}+{}\\
{}+\dfrac{2\alpha(2\alpha+2n_1+1)(2\alpha-2n_1-1)}{2\alpha+1}\cdot
\frac{b(b-1)(b-a)B}{(z-b)^{2\alpha+1}}+{}\\
{}+O\left((z-b)^{-2\alpha}\right), \qquad z\to b,
\end{multline*}
%
has at the same point a pole of order $\alpha'\leqslant 2n_1+1$.

In particular
\begin{multline*}
I_1=-3n_1(n_1-1)(n_1+1)(n_1+2)b(b-1)(b-a)\times{}\\
\times\left(\frac{b(b-1)(b-a)}{(z-b)^4}-
\dfrac{3b^2-2(a+1)b+a}{(z-b)^3}\right)+O\left((z-b)^{-2}\right).
\end{multline*}

Let now $n_1=1$. Then the function  $I_1$ has a pole of second
order at the point $z=b$ while $I_0$  has also a pole of second
order at the same point. Hence, there exist a linear combination
$\wtI_1=I_1+\wtc I_0$ with a pole of first order in $s=b$.
Let us find $b$ from condition
\begin{equation} \label{eq:b}
\Res_{z=b}\wtI_1=0.
\end{equation}

Let us take functions $\wtI_j,$
\begin{equation}
\quad \wtI_j=\Ll(\wtI_{j-1}),\quad  j>1.
\end{equation}
%
Firstly this functions do not have poles at the point $z=b$.
%
Secondly, for any $j$ the order of poles of functions $\wtI_j(z)$
in singular points $z=z_k$ ($ z_1=0, z_2=1, z_3=a$) of the
Fuchsian equation does not exceed corresponding
characteristic~$m_k$ (see for example \cite{Sm2002}).
%

Hence, for any $j$ the dimension of the linear span of rational
functions $1,\wtI_1,\ldots,\wtI_j$ does not exceed $N$ and
therefore there exists a number $g$,
\begin{equation}
\max_{0\leqslant i\leqslant 3} m_i \leqslant g\leqslant N-1,
\label{eq:g}
\end{equation}
%
such that the equality \eqref{eq:nov1} is fulfilled.
%
We remind that functions $ \wtI_j$ are linear combinations of
functions $I_i$ of lesser orders (Remark~\ref{rem1}).
\end{proof}

The condition \eqref{eq:b} can be found in explicit form. Direct
calculation give us next polynomial equation of sixth order on~$b$
(position of false singular point):
\begin{equation} \label{eq:b1}
k_0^2 b^6-2k_0^2(a+1)b^5+h_4 b^4+2ah_3 b^3+ ah_2
b^2+2k_1^2a^2(a+1)b-k_1^2a^3=0,
\end{equation}
%
where $k_j=m_j-1/2$
\begin{gather*}
h_4=(k_0^2-k_3^2)a^2+(4k_0^2+k_2^2+k_3^2-k_1^2)a+k_0^2-k_2^2,\\
h_3=(k_1^2+k_3^2-k_0^2-k_2^2)a+k_1^2+k_2^2-k_0^2-k_3^2,\\
h_2=(k_2^2-k_1^2)a^2+(k_0^2-4k_1^2-k_2^2-k_3^2)a+k_3^2-k_1^2.
\end{gather*}

\begin{rem}
Any value  $b$  satisfying the sixth-order polynomial equation
\eqref{eq:b1} gives a `finite-gap' Fuchsian equation and a
finite-gap elliptic potential of Shr\"odinger operator. Hence, in
general cases for any set of $m_i\in\Bz_{\geqslant 0}$ there exist
six different finite-gap elliptic potentials ($s=1,\ldots,6$):
\begin{equation}
u(x)=m_0(m_0+1)\wp(x)+\sum_{i=1}^3m_i (m_i+1)\wp(x-\omega_i) +
2\wp(x-\delta_s)+2\wp(x+\delta_s).
\end{equation}
For example, if
\begin{equation*}
m_0=m_1=m_2=m_3=0, \quad n_1=1
\end{equation*}
then
\begin{equation*}
u(x)= 2\wp(x-\delta_s)+2\wp(x+\delta_s),
\end{equation*}
where
\begin{equation*}
\delta_1=\frac{\omega}{2},\quad
\delta_2=\frac{\omega}{2}+\omega',\quad
\delta_3=\frac{\omega'}{2},\quad
\delta_4=\omega+\frac{\omega'}{2},\quad
\delta_5=\frac{\omega+\omega'}{2}, \quad
\delta_6=\frac{\omega-\omega'}{2}.
\end{equation*}
\end{rem}

\begin{rem}
Another proof of theorem~\ref{thm3} is based on
theorem~\ref{thm1}. The condition of absence of logarithmic
singularity of general solution of equation \eqref{eq:main},
\eqref{eq:P:five}--\eqref{eq:N:five} imposes a restriction on the
position of the point~$b$. This restriction is equivalent to
equation \eqref{eq:b1}.
\end{rem}

\begin{rem}\label{rem3}
There are special cases which show that conditions $M=1$, $n_1=1$
are sufficient but not necessary. In particular if
\begin{equation*}
m_0=m_1=m_2=m_3=0, \quad n_1=g,\quad  \delta_1=\frac{\omega_j}2
\end{equation*}
%
the potential
\begin{equation*}
u(x)=g(g+1)\wp(x-\omega_j/2)+g(g+1)\wp(x+\omega_j/2)
\end{equation*}
is $g$-gap potential because it represents a $g$-gap Lam\'e
potential with changed period of lattice
\begin{equation*}
u(x)=g(g+1)\wp(x)+g(g+1)\wp(x+\omega_j)
\end{equation*}
which shifted on one fourth of new period.
With the help of formula \eqref{eq:multi} and with the shift of
argument one can construct from Treibich-Verdier potentials
\eqref{pot:tv} more complicated even finite-gap potentials. Hence,
there it exists more complicated (but more special) `finite-gap'
Fuchsian equations.
\end{rem}

In order to find finite-gap solutions of Fuchsian equation
\eqref{eq:main}, \eqref{eq:P:six}--\eqref{eq:N:six} let us use the
results of the theory of finite-gap elliptic potentials for the
Schr\"odinger operator and consider the equation
\begin{equation}
\frac{d^3\Psi}{dz^3}+3P(z)\frac{d^2\Psi}{dz^2}
+\left(P'(z)+4Q(z)+2P^2(z)\right)\frac{d\Psi}{d z}
+\left(2Q'(z)+4P(z)Q(z)\right)\Psi=0,\label{eq:3}
\end{equation}
%
solutions of which are the products of any two solutions of Heun's
equation \eqref{eq:main}.

\begin{thm} \label{thm:six:Psi}
If the equation \eqref{eq:main},
\eqref{eq:P:six}--\eqref{eq:N:six} is finite-gap then the equation
\eqref{eq:3}, \eqref{eq:P:six}--\eqref{eq:N:six} with nonnegative
integer characteristics $m_i,n_j$ has as its solution the function
$\Psi_{g,N}(\lambda,z)$, which is a polynomial in $\lambda$ of the
degree $g$ and in $z$ of the degree $N$ \eqref{eq:N:six}
\begin{equation} \label{eq:six:Psi}
\begin{split}
\Psi_{g,N}(\lambda,z)&=a_0(\lambda) z^N+a_1(\lambda) z^{N-1}
+\ldots+a_N(\lambda)={}\\
&{}=\wta_0(z)\lambda^g +\wta_1(z)\lambda^{g-1}+\ldots+\wta_g(z).
\end{split}
\end{equation}
%
The leading coefficient of this function considered as a
polynomial in $\lambda$ is equal to
\begin{equation}
\wta_0(z)=z^{m_1}(z-1)^{m_2}(z-a)^{m_3}\prod_{k=1}^M(z-b_k)^{2n_k}.
\label{six:a0}
\end{equation}
\end{thm}

\begin{proof}
The product $\whPsi(x,E)$ of eigenfunctions of the Schr\"odinger
operator with the potential $u(x)$ \eqref{pot:gen} is an elliptic
meromorphic function in the variable $x$, because $u(x)$ is a
Picard potential (Proposition~\ref{utv:picard}).
As a function of the variable $x$ the function $\whPsi(x,E)$ has
poles of multiplicity $2m_j$ at the points $x=\omega_j$
($\omega_0\equiv 0$) and $2n_k$ at the points $x=\pm\delta_k$.

From equation \eqref{kdv:prod:psi} and from evenness $u(x)$
\eqref{pot:gen} it follow what $\whPsi(x,E)$ is rational function
in $\wp(x)$ and a polynomial of degree $g$ in the spectral
parameter $E$. Hence, the function
\begin{equation}
\Psi_{g,N}(\lambda,z)=\const\cdot\whPsi(x,E)
\prod_{j=1}^3(\wp(x)-e_j)^{m_i}\prod_{k=1}^M(\wp(x)-\wp(\delta_k))^{2n_k},
\label{six:Psi:def}
\end{equation}
%
where $\lambda$ and $z$ are related with $E$ and $x$ by equalities
\eqref{eq:x2z:six}, is a polynomial in $\lambda$ of degree $g$ and
in $z$ of degree $N$ (i.e. it is a rational function in the
variable $z$ with the unique pole of order $N$ at the point
$z=\infty$). The constant in the equality \eqref{six:Psi:def} is
chosen such that the leading coefficient of
$\Psi_{g,N}(\lambda,z)$, considered as a polynomial in $\lambda$,
is equal \eqref{six:a0}.
\end{proof}

\begin{cor} Coefficients $\wta_j(z)$ of the polynomial $\Psi_{g,N}(\lambda,z)$
have the form:
\begin{equation}
\wta_j(z)=\wta_0(z)\whI_{j-1}, \qquad j=1,\ldots,g,
\end{equation}
where $\wtI_{j}$ is a linear combination of the rational functions
$I_j,\ldots,I_0,1$ having poles in the singularities of Fuchsian
equation.%
\end{cor}

\begin{proof}[Proof\nodot] follows from the
equalities \eqref{kdv:gamma}, \eqref{six:Psi:def} and from the
change of variable \eqref{eq:x2z:six}.
\end{proof}

\begin{cor}
If the equation \eqref{eq:3} does not have solution in polynomial
form \eqref{eq:six:Psi}, then the Fuchsian equation
\eqref{eq:main}, \eqref{eq:P:six}--\eqref{eq:N:six} is not
finite-gap.
\end{cor}

Knowing the product of the solutions of Heun's equation it is not
difficult to find the solutions themselves (see, for instance,
\cite[\S 19.53, \S 23.7, \S23.71]{WW}).

\begin{thm} Finite-gap solutions of Fuchsian equation \eqref{eq:main},
\eqref{eq:P:six}--\eqref{eq:N:six} with $m_i,n_k\in\Bz_{\geqslant 0}$
have the form
\begin{multline} \label{eq:Y}
Y_{1,2}(\mathbf{m},\mathbf{n};\lambda;z)={}\\
{}=\sqrt{\Psi_{g,N}(\lambda,z)}\exp\left(\pm\frac{\ii\nu(\lambda)}2
\int\frac{z^{m_1}(z-1)^{m_2}(z-a)^{m_3}\prod_{k=1}^M(z-b_k)^{2n_k}\,d
z}{\Psi_{g,N}(\lambda,z) \sqrt{z(z-1)(z-a)}}\right).
\end{multline}
Here $\ii^2=-1$,
\begin{equation}
\vG:\quad \nu^2=\prod_{j=1}^{2g+1}(\lambda-\lambda_j), \qquad
\lambda_j=\lambda(E_j), \label{curve:hyp}
\end{equation}
$E_j$ are the gap edges of the finite-gap elliptic potential
$u(x)$ \eqref{pot:gen}.
\end{thm}

\begin{proof}
From the Liouville formula it follows that the Wronskian of two
linearly independent solutions of the linear homogeneous
differential equation,
\begin{equation*}
y''+P(z)y'+Q(z)y=0,
\end{equation*}
has the following dependence from $z$:
\begin{equation*}
W[y_1(\lambda,z),y_2(\lambda,z)]=W_0(\lambda)
\exp\left\{-\int_{z_0}^z P(t)\,dt\right\},
\end{equation*}
where $W_0(\lambda)$ is Wronskian's value at $z_0$. Hence, the
Wronskian of two linearly independent solutions of Fuchsian
equation \eqref{eq:main}, \eqref{eq:P:six}--\eqref{eq:N:six} is
equal to
\begin{equation}\label{eq:wron}
W[y_1(\lambda,z),y_2(\lambda,z)]=-\ii\nu(\lambda)\cdot
\frac{z^{m_1}(z-1)^{m_2}(z-a)^{m_3}\prod_{k=1}^M(z-b_k)^{2n_k}}{\sqrt{z(z-1)(z-a)}},
\end{equation}
where
\begin{equation*}
\nu(\lambda)= \ii W_0(\lambda)\cdot\frac{\sqrt{z_0(z_0-1)(z_0-a)}
}{z_0^{m_1}(z_0-1)^{m_2}(z_0-a)^{m_3}\prod_{k=1}^M(z_0-b_k)^{2n_k}}.
\end{equation*}
If we now divide the Wronskian of two solutions \eqref{eq:wron} by
their product,
\begin{equation}\label{eq:prod}
y_1(\lambda,z)\cdot y_2(\lambda,z)=\Psi_{g,N}(\lambda,z),
\end{equation}
we obtain a simple differential equation of the first order:
\begin{equation}\label{eq:frac:dif}
\frac{y_2'}{y_2}-\frac{y_1'}{y_1}= -\ii\nu(\lambda)\cdot
\frac{z^{m_1}(z-1)^{m_2}(z-a)^{m_3}\prod_{k=1}^M(z-b_k)^{2n_k}}%
{\Psi_{g,N}(\lambda,z)\sqrt{z(z-1)(z-a)}}.
\end{equation}
The equation \eqref{eq:frac:dif} can be easily integrated:
\begin{equation}\label{eq:frac}
\frac{y_2(\lambda,z)}{y_1(\lambda,z)}
=C(\lambda)\cdot\exp\left[-\ii\nu(\lambda)\int_{z_1}^z
\frac{t^{m_1}(t-1)^{m_2}(t-a)^{m_3}
\prod_{k=1}^M(t-b_k)^{2n_k}\,dt}{\Psi_{g,N}(\lambda,t)
\sqrt{t(t-1)(t-a)}}\right],
\end{equation}
where $C(\lambda)=y_2(\lambda,z_1)/y_1(\lambda,z_1)$.

If we now consider solutions
$Y_{1,2}(\mathbf{m},\mathbf{n};\lambda,z)$ with the same Wronskian
\eqref{eq:wron} and product \eqref{eq:prod}
\begin{equation}
Y_1(\mathbf{m},\mathbf{n};\lambda,z)=\sqrt{C(\lambda)}\cdot
y_1(\lambda,z),\qquad Y_2(\mathbf{m},\mathbf{n};\lambda,z)=
\frac{y_2(\lambda,z)}{\sqrt{C(\lambda)}},
\end{equation}
then from \eqref{eq:prod} and \eqref{eq:frac} we get \eqref{eq:Y}.

Substituting the ansatz \eqref{eq:Y} into the Fuchsian equation
\eqref{eq:main}, \eqref{eq:P:six}--\eqref{eq:N:six} we get
\begin{equation}
\nu^2(\lambda)=\frac{2\Psi\Psi''-(\Psi')^2+2P(z)\Psi\Psi'
+4Q(z)\Psi^2}{z^{2m_1-1}(z-1)^{2m_2-1}(z-a)^{2m_3-1}
\prod_{k=1}^M(z-b_k)^{4n_k}}, \label{eq:six:nu}
\end{equation}
where
\begin{equation*}
\Psi=\Psi_{g,N}(\lambda,z),\qquad \Psi'= \frac{d\Psi}{d z},\qquad
\Psi''= \frac{d^2\Psi}{d z^2}.
\end{equation*}
From \eqref{eq:six:Psi}, \eqref{eq:six:nu} it follows that
$\nu^2(\lambda)$ is a polynomial in $\lambda$ of the degree $2g+1$
with the leading coefficient equal to~$1$.

It is not difficult to show that under the change
\eqref{eq:x2z:six} the solutions
$Y_{1,2}(\mathbf{m},\mathbf{n};\lambda,z)$ of Fuchsian equation
\eqref{eq:main}, \eqref{eq:P:six}--\eqref{eq:N:six} turn into
Floquet solutions of the equation \eqref{eq:shr}, \eqref{pot:gen}.
Therefore, zeros $\lambda_j$ ($j=1,\ldots,2g+1$) of the polynomial
$\nu(\lambda)$ or, which is the same, zeros of the Wronskian of
the solutions $Y_{1,2}(\mathbf{m},\mathbf{n};\lambda,z)$
correspond to zeros of the Wronskian of the Floquet solutions of
equation \eqref{eq:shr}, \eqref{pot:gen}, i.e. they correspond to
the gap edges $E_j$ ($j=1,\ldots,2g+1$) of spectrum of the
potential \eqref{pot:gen}. Hence, the hyperelliptic curve $\vG$
\eqref{curve:hyp} is isomorphic to the spectral curve $\wtG$
\begin{equation} \label{curve:hyp1}
w^2=\prod_{j=1}^{2g+1}(E-E_j)
\end{equation}
of the finite-gap elliptic potential $u(x)$ \eqref{pot:gen}.
\end{proof}

\begin{rem}
Knowing the product of the eigenfunctions of the Shr\"odinger
operator $\whPsi(x,E)$ \eqref{kdv:prod:psi} it is possible to
write the formulae for these eigenfunctions
\begin{equation}
\psi_{1,2}(x,E)=\sqrt{\whPsi(x,E)}\exp\left(\mp
w(E)\int\dfrac{dx}{\whPsi(x,E)}\right), \qquad
w(E)=\dfrac{1}{2}W[\psi_1,\psi_2]
\end{equation}
%
and for equation of spectral curve \eqref{curve:hyp1}
\begin{equation}
w^2(E)=(u+E)\whPsi^2+\dfrac{1}{4}\whPsi^2_x-\dfrac{1}{2}\whPsi_{xx}\whPsi.
\end{equation}
\end{rem}

\section*{Concluding remarks}

There is one important difference between finite-gap Heun equation
\cite{Sm2002} and finite-gap Fuchsian equation. In case of the
Heun equation the spectral curve and its algebraic genus $g$ are
completely determined by characteristics $m_i$:
\begin{enumerate}
\item in the case of even $N=\sum m_i$
\begin{equation*}
g=\max\left\{\max_{0\leqslant i\leqslant3}m_i, \frac N2
-\min_{0\leqslant i\leqslant3}m_i\right\};
\end{equation*}
\item in the case of odd $N$
\begin{equation*}
g=\max\left\{\max_{0\leqslant i\leqslant3}m_i, \frac
{N+1}2\right\}.
\end{equation*}
\end{enumerate}

On the other hand in case of the Fuchsian equation the spectral
curve and its algebraic genus $g$ depend not only on
characteristics $m_i$, $n_k$ but also on positions of singular
points $b_k$.
For example, the potential
\begin{equation*}
u(x)=2\wp(x)+2\wp(x+\delta)+2\wp(x-\delta)
\end{equation*}
%
for $\wp(2\delta)=-2\wp(\delta)$ is two-gap, but for
$\wp(2\delta)=\wp(\delta)$ is one-gap potential. Therefore,
knowing only characteristics $m_i$, $n_k$ we can only estimate
(see \eqref{eq:g}) algebraic genus of spectral curve
\eqref{curve:hyp}.
%
Examples of simple `finite-gap' solutions of Fuchsian equation and
of finite-gap elliptic potentials that are not Lam\'e or
Treibich-Verdier potentials, are discussed in Appendix.

It is easy to see that equation \eqref{eq:shr}, \eqref{pot:gen} is
invariant with respect to transformation $m_i \to -m_i-1$, $n_k \to
-n_k-1$.  Therefore it is not difficult to transform `finite-gap'
solutions of equation \eqref{eq:main}, \eqref{eq:P:six}--\eqref{eq:N:six}
with non-negative characteristics $m_i$, $n_k$ into solutions with
negative characteristics. Corresponding transformations for Heun equation
one can find, for example, in \cite{ODE, VTF, HDE,
SL, Sm2002}.

\appendix
\section{Simplest `finite-gap' solutions}

At the end of this paper we would like to list equations for $b$,
polynomials $\Psi(\lambda,z)$ and canonical equations
\eqref{curve:hyp} of the hyperelliptic curves
$\vG=\{(\nu,\lambda)\}$ for some simplest `finite-gap' solutions
\eqref{eq:Y} of Fuchsian equation with five singular points
\eqref{eq:main}, \eqref{eq:P:five}--\eqref{eq:N:five} with
characteristics $n_1=1$, $\quad m_i\in\Bz_{\geqslant 0}$ $(i=0,1,2,3)$.
%
There will be given also finite-gap elliptic potentials
$\wtu(x)=u(x)+\const$ and their spectral curves
\eqref{curve:hyp1}. Potentials $\wtu(x)$ are normalized by
condition
\begin{equation*}
\sum_{j=1}^{2g+1} \wtE_j=0.
\end{equation*}
Our examples are indexed by
the characteristics $(m_0,m_1,m_2,m_2)$.

\begin{description}
\item[(0,0,0,0)]
\begin{equation}
(b^2-a)(b^2-2b+a)(b^2-2ab+a)=0.
\end{equation}
$\mathbf{b^2-a=0}$:
\begin{equation*}
\Psi_{1,2}(\lambda,z)=(z-b)^2\lambda+(3+3a-4b)z^2-2(5b+5ab-8a)z+a(3+3a-4b),
\end{equation*}
\begin{equation*}
\nu^2=(\lambda-4b+3a+3)(\lambda^2+7(1+a-2b)\lambda+2(6a^2+36a+6-25ab-25b)).
\end{equation*}
%
The potential $ \wtu(x) $ is 
one-gap potential with changed period of lattice.
\begin{align*}
\wtu(x)&=2\wp(x-\delta)+2\wp(x+\delta)-2e_1, \qquad\wp(2\delta)=e_1,\\
w^2&=(E+2e_1)(E+e_1-2\wp(\delta))(E-3e_1+2\wp(\delta)).
\end{align*}
$\mathbf{b^2-2b+a=0}$:
\begin{equation*}
\Psi_{1,2}(\lambda,z)=(z-b)^2\lambda+(3a-4b)z^2-2(6a+5ab-12b)z-3a^2+12a-24b+14ab,\\
\end{equation*}
\begin{equation*}
\nu^2=(\lambda-4b+3a)(\lambda^2+(4+7a-14b)\lambda+2(6a^2-16a-25ab+36b)),\\
\end{equation*}
%
The potential $ \wtu(x) $ is 
one-gap potential with changed period of lattice.
\begin{align*}
\wtu(x)&=2\wp(x-\delta)+2\wp(x+\delta)-2e_2, \qquad\wp(2\delta)=e_2,\\
w^2&=(E+2e_2)(E+e_2-2\wp(\delta))(E-3e_2+2\wp(\delta)).
\end{align*}
$\mathbf{b^2-2ab+a=0}$:
\begin{equation*}
\Psi_{1,2}(\lambda,z)=(z-b)^2\lambda+(3-4b)z^2-2(6a+5b-12ab)z-3a+12a^2-24a^2b+14ab,\\
\end{equation*}
\begin{equation*}
\nu^2=(\lambda-4b+3)(\lambda^2+(4a+7-14b)\lambda+2(6-16a-25b+36ab)),\\
\end{equation*}
%
The potential $ \wtu(x) $ is 
one-gap potential with changed period of lattice.
\begin{align*}
\wtu(x)&=2\wp(x-\delta)+2\wp(x+\delta)-2e_3, \qquad\wp(2\delta)=e_3,\\
w^2&=(E+2e_3)(E+e_3-2\wp(\delta))(E-3e_3+2\wp(\delta)).
\end{align*}
\end{description}

\begin{description}
\item[(1,0,0,0)]
\begin{equation}
(3b^2-2(a+1)b+a)(3b^4-4(a+1)b^3+6ab^2-a^2)=0.
\end{equation}
$\mathbf{3b^2-2(a+1)b+a=0}$:
\begin{multline*}
\Psi_{2,3}(\lambda,z)=(z-b)^2\lambda^2+\left(z^3+3(1+a-4b)z^2
-(7a-8ab-8ab)z
+\vphantom{\dfrac{1}{3}}\right.\\
{}+\left.\dfrac{1}{3}(5a^2+5a-10a^2b-10b-2ab)\right)\lambda-\\
{}-\dfrac{4}{3}(a^3-4a^2+a-2a^3b+3a^2b+3ab-2b) ,
\end{multline*}
\begin{multline*}
\nu^2=\lambda(\lambda^4+10(a+1-3b)\lambda^3+33(a^2-a+1)\lambda^2+\\
{}+3(12a^3+a^2+a+12-38a^2b+38ab-38b)\lambda-\\
{}-12(a^3-4a^2+a-2a^3b+3a^2b+3ab-2b)).
\end{multline*}
%
The potential $ \wtu(x) $ is an isospectral deformation of two-gap
Lam\'e potential.
\begin{align*}
\wtu(x)&=2\wp(x)+2\wp(x-\delta)+2\wp(x+\delta),
\qquad\wp(2\delta)=-2\wp(\delta),
\quad \wp^2(\delta)=\frac{g_2}{12},\\
w^2&=\left(E^3-\dfrac{9g_2}{4}E+\dfrac{27g_3}{4}\right)
(E-6\wp(\delta))(E+6\wp(\delta)).
\end{align*}
$\mathbf{3b^4-4(a+1)b^3+6ab^2-a^2=0}$:
\begin{multline*}
\Psi_{1,3}(\lambda,z)=(z-b)^2\lambda+z^3+3(1+a-4b)z^2+\\
{}+(2a-10ab-10b+27b^2)z+2ab+3ab^2+3b^2-12b^3,
\end{multline*}
\begin{multline*}
\nu^2=\lambda^3+10(1+a-3b)\lambda^2+3(11a^2+19a+90b^2-60ab-60b+11)\lambda+\\
{}+36a^3+73a^2+73a+36-254b-376ab-254a^2b+630b^2+630ab^2-630b^3.
\end{multline*}
%
The potential $ \wtu(x) $ is one-gap potential with changed period
of lattice.
\begin{align*}
\wtu(x)&=2\wp(x)+2\wp(x-\delta)+2\wp(x+\delta)-4\wp(\delta),
\qquad\wp(2\delta)=\wp(\delta),\\
w^2&=E^3+\left(\dfrac{9g_2}{4}-30\wp^2(\delta)\right)E
-70\wp^3(\delta)+\dfrac{21g_2}{2}\wp(\delta)+\dfrac{27g_3}{4}.
\end{align*}
\end{description}

\begin{description}
\item[(1,1,0,0)]
\begin{equation}
(b^2-a)(3b^2-2(a+2)b+3a)(3b^2-2(2a+1)b+3a)=0.
\end{equation}
$\mathbf{b^2-a=0}$:
\begin{multline*}
\Psi_{1,4}(\lambda,z)=z(z-b)^2\lambda+z^4+8(1+a-2b)z^3+\\
{}+2(19a-10ab-10b)z^2+8a(1+a-2b)z+a^2,
\end{multline*}
\begin{multline*}
\nu^2=\lambda^3+(25a+25-42b)\lambda^2+8(26a^2+120a+26-85b-85ab)\lambda+\\
{}+144(4a^3+41a^2+41a+4-19a^2b-52ab-19b).
\end{multline*}
%
The potential $ \wtu(x) $ is one-gap potential with changed period
of lattice.
\begin{align*}
\wtu(x)&=2\wp(x)+2\wp(x-\omega_1)+2\wp(x-\delta)+2\wp(x+\delta)
-4\wp(\delta)-2e_1,
\qquad\wp(2\delta)=e_1,\\
w^2&=E^3-(75e_1^2+60e_1\wp(\delta)-11g_2)E-2e_1(137e_1^2-15g_2)
-28(15e_1^2-g_2)\wp(\delta).
\end{align*}
$\mathbf{3b^2-2(a+2)b+3a=0}$:
\begin{multline*}
\Psi_{2,4}(\lambda,z)=z(z-b)^2\lambda^2+\left(z^4+(13a+8-20b)z^3
-\dfrac{4}{3}(30a+ab-28b)z^2\right.+\\
{}+\left.\dfrac{1}{9}(3a^2+168a-2a^2b+64ab-224b)z
-\dfrac{a}{3}(3a-2ab-4b)\right)\lambda+\\
{}+(5a-4b)z^4+\dfrac{8}{3}(15a^2-9a-26ab+20b)z^3-\\
{}-\dfrac{2}{9}(660a^2-480a-35a^2b-776ab+640b)z^2+\\
{}+\dfrac{4}{27}(15a^3+492a^2-480a-10a^3b+39a^2b-696ab+640b)z+\\
{}+\dfrac{a}{9}(15a^2+120a-10a^2b+44ab-160b),
\end{multline*}
\begin{multline*}
\nu^2=(\lambda+5a-4b)\left(\lambda^4+(30a+25-46b)\lambda^3+
(333a^2-199a+208-506ab+188b)\lambda^2\vphantom{\dfrac{1}{2}}\right.+\\
{}+\dfrac{4}{9}(3636a^3-7212a^2+3612a+1296-5449a^2b+
9065ab-4912b)\lambda+\\
{}+\dfrac{8}{27}(9720a^4-32019a^3+26751a^2-3912a-\\
{}-\left.\vphantom{\dfrac{1}{2}}14555a^3b+44451a^2b-40836ab+10400b)\right).
\end{multline*}
%
The potential  $\wtu(x)$ is an isospectral deformation of two-gap
$4$-elliptic Treibich-Verdier potential with additional pole at
the point $x=\omega_2$ (see, for example,
\cite{Smr89a,Smr94,BE89b,BE94}).
\begin{align*}
\wtu(x)&=2\wp(x)+2\wp(x-\omega_1)+2\wp(x-\delta)+2\wp(x+\delta)
-2e_2,
\qquad\wp(2\delta)=\dfrac{e_1+5e_2}{3},\\
w^2&=(E-6e_2)(E^2+2E(e_2-e_3)-3(13e_2^2+2e_2e_3-3e_3^2)\times\\
&\quad
{}\times(E-3(e_1-e_2)-6\wp(\delta))(E+e_1-e_2+6\wp(\delta)).
\end{align*}
$\mathbf{3b^2-2(2a+1)b+3a=0}$:
\begin{multline*}
\Psi_{2,4}(\lambda,z)=z(z-b)^2\lambda^2+\left(z^4+(13+8a-20b)z^3
-\dfrac{4}{3}(30a+b-28ab)z^2\right.+\\
{}+\left.\dfrac{1}{9}(3a+168a^2-2b+64ab-224a^2b)z
-\dfrac{a}{3}(3a-2b-4ab)\right)\lambda+\\
{}+(5a-4b)z^4+\dfrac{8}{3}(15-9a-26b+20ab)z^3-\\
{}-\dfrac{2}{9}(660a-480a^2-35b-776ab+640a^2b)z^2+\\
{}+\dfrac{4}{27}(15a+492a^2-480a^3-10b+39ab-696a^2b+640a^3b)z+\\
{}+\dfrac{a}{9}(15a+120a^2-10b+44ab-160a^2b),
\end{multline*}
\begin{multline*}
\nu^2=(\lambda+5-4b)\left(\lambda^4+(30+25a-46b)\lambda^3+
(333-199a+208a^2-506b+188ab)\lambda^2\vphantom{\dfrac{1}{2}}\right.+\\
{}+\dfrac{4}{9}(3636-7212a+3612a^2+1296a^3-5449b+
9065ab-4912a^2b)\lambda+\\
{}+\dfrac{8}{27}(9720-32019a+26751a^2-3912a^3-\\
{}-\left.\vphantom{\dfrac{1}{2}}14555b+44451ab-40836a^2b+10400a^3b)\right).
\end{multline*}
%
The potential  $\wtu(x)$ is an isospectral deformation of two-gap
$4$-elliptic Treibich-Verdier potential
\cite{Smr89a,Smr94,BE89b,BE94} with additional pole at the point
$x=\omega_3$
\begin{align*}
\wtu(x)&=2\wp(x)+2\wp(x-\omega_1)+2\wp(x-\delta)+2\wp(x+\delta)
-2e_3,
\qquad\wp(2\delta)=\dfrac{e_1+5e_3}{3},\\
w^2&=(E-6e_3)(E^2-2E(e_2-e_3)-3(13e_3^2+2e_2e_3-3e_2^2)\times\\
&\quad
{}\times(E-3(e_1-e_3)-6\wp(\delta))(E+e_1-e_3+6\wp(\delta)).
\end{align*}
\end{description}

\begin{description}
\item[(2,0,0,0)]
\begin{multline}
25b^6-50(a+1)b^5+(24a^2+101a+24)b^4-\\
{}-48a(a+1)b^3+19a^2b^2+2a^2(a+1)b-a^3=0.
\end{multline}
\begin{multline*}
\Psi_{2,4}(\lambda,z)=(z-b)^2\lambda^2+(3z^3+(3a+3-26b)z^2+\\
{}+(2a-10ab-10b+49b^2)z+b(2a+3ab+3b-22b^2))\lambda+9z^4-48bz^3-\\
{}-(9a-6ab-6b-115b^2)z^2+2b(17a-10ab-10b-61b^2)z+\\
{}+3a^2-59ab^2+42ab^3+42b^3+21b^4,
\end{multline*}
\begin{multline*}
\nu^2=\lambda^5+10(a+1-5b)\lambda^4+(33a^2+17a+33-260ab-260b+790b^2)\lambda^3+\\
{}+(36a^3-135a^2-135a+36-234a^2b+\\
{}+624ab+1170ab^2-234b+1170b^2-4330b^3)\lambda^2-\\
{}-(216a^3-189a^2+216a-144b+108ab+108a^2b-144a^3b-\\
-828b^2-924ab^2-828a^2b^2+4300b^3+4300ab^3-11885b^4)\lambda-\\
{}-4(27a^2+27a^3+108ab-351a^2b+108a^3b-36b^2-405ab^2-405a^2b^2-\\
-36a^3b^2+528b^3+3008ab^3+582a^2b^3-2870b^4-2870ab^4+4099b^5).
\end{multline*}
%
The potential  $\wtu(x)$ is a finite-gap elliptic potential which
can not be transformed into Lam\'e or Treibich-Verdier potential
by shifting of argument or by transformation of lattice of periods
\eqref{eq:multi}. It seems that this potential first appeared in
\cite{Smr89a,Smr94}.
\begin{equation*}
\wtu(x)=6\wp(x)+2\wp(x-\delta)+2\wp(x+\delta)-4\wp(\delta), \qquad
\wp'(2\delta)=-3\wp'(\delta),
\end{equation*}
\begin{multline*}
w^2=E^5-\left(210\wp^2(\delta)-\dfrac{49g_2}{4}\right)E^3
+\left(630\wp^3(\delta)-\dfrac{189g_2}{2}\wp(\delta)-\dfrac{225g_3}{4}\right)E^2+\\
{}+\left(12285\wp^4(\delta)-\dfrac{3213g_2}{2}\wp^2(\delta)
-297g_3\wp(\delta)+\dfrac{621g_2^2}{16}\right)E-\\
{}-59454\wp^5(\delta)+11421g_2\wp^3(\delta) +3240g_3\wp^2(\delta)
-\dfrac{4239g_2^2}{8}\wp(\delta)-\dfrac{1107g_2g_3}{4}.
\end{multline*}
\end{description}

\end{document}